\input amstex
\magnification=\magstep 1

\TagsOnRight
\documentstyle{amsppt}
\overfullrule=0pt \vsize=550pt

\define\wre{with respect to}
\define\sqe{sequence}

\define\bsk{\bigskip}

\define\ineq{inequality}
\define\ineqs{inequalities}

\define\cd{conditional distribution}

\define\e{\varepsilon}

\define\Lw{\Cal L}
\define\si{\sigma}
\define\Ga{\Gamma}

        \define\jd{joint distribution}

\define\rv{random variable}

\define\de{\delta}

\define\LSI{logarithmic Sobolev \ineq}
\define\LSIs{logarithmic Sobolev inequalities}

\define\sumk{\sum_{k=1}^n}
\define\R{\Bbb R}
\define\qute{\quad\text}

\define\bV{\bigm\Vert}

\define\pms{probability measure}

\document

\topmatter

\title
{An inequality for relative entropy and logarithmic Sobolev inequalities in Euclidean spaces}
\endtitle

 \leftheadtext {An inequality for relative entropy}
\rightheadtext {logarithmic Sobolev inequalities}

\author Katalin Marton \endauthor
\affil Alfr\'ed R\'enyi  Institute of Mathematics,
Hungarian Academy of Sciences
\endaffil

\date{June 18, 2012}
\enddate

\address H-1364 POB 127, Budapest, tel. 36 (1) 4838300, Hungary
\endaddress
\email marton\@renyi.hu
\endemail
\footnote{This work was supported in part by the grants OTKA
T-32323 and K-76088 of the
Hungarian Academy of Sciences}
\footnote{ AMS 2000 subject  classifications: 52A40, 60K35, 82B20, 82C22.}
\footnote{Key words and phrases: Relative entropy, Wasserstein distance,
Fokker-Planck equation,
gradient flow,
non-compact spin system, Gibbs sampler,
weakly dependent random variables,
logarithmic Sobolev inequality,
transportation-cost inequality.}
\footnote{Subject classification:
82C22, 
60J05, 
35Q84, 
60J25, 
82B21.} 

\abstract
For a class of density functions $q(x)$ on $\Bbb R^n$
we prove
an inequality between relative entropy and the
weighted sum of  conditional
relative entropies of the following form:
$$
\align
&D(p||q)\leq\\
& Const.
\sum_{i=1}^n
\rho_i\cdot
D(p_i(\cdot|Y_1,\dots, Y_{i-1},Y_{i+1},\dots, Y_n)
||
  Q_i(\cdot|Y_1,\dots, Y_{i-1},Y_{i+1},\dots, Y_n))
\endalign
$$
for any density function
$p(x)$ on $\Bbb R^n$,
where $p_i(\cdot|y_1,\dots, y_{i-1},y_{i+1},\dots, y_n)$ and
\newline
$Q_i(\cdot|x_1,\dots, x_{i-1},x_{i+1},\dots, x_n)$ denote
the local specifications of $p$ resp. $q$, and
$\rho_i$ is
the logarithmic Sobolev constant of
$Q_i(\cdot|x_1,\dots, x_{i-1},x_{i+1},\dots, x_n)$.
Thereby we derive a logarithmic Sobolev inequality
for a weighted Gibbs sampler governed by the
local specifications of $q$. Moreover, the above
 inequality implies a classical
logarithmic Sobolev inequality for $q$,
as defined for Gaussian distribution by L. Gross.
This strengthens
a result by F. Otto and M. Reznikoff.
The proof is based on  ideas
developed by F. Otto and C. Villani in their paper on
 the connection between Talagrand's
transportation-cost inequality and logarithmic Sobolev inequality.

\endabstract
\endtopmatter

\beginsection 1. Introduction.

The motivation for this paper was to prove \LSIs\ on product spaces, under
possibly general conditions.
\bsk

First we define some basic concepts:
\bsk

\definition{Definition}
For \pms s $p$ and $q$ on $\R^m$ ($m\ge 1$ integer),
we denote by $D(p\Vert q)$ the relative entropy  of  $p$ \wre\   $q$:
$$
D(p\Vert q) = \int_{\R^m} \log \frac {dp(u)} {dp(u)}dp(u) \qute{if}\quad p<<q,
\tag 1.1
$$
and $\infty$ otherwise.
We always have in mind \pms s absolutely continuous \wre\ the Lebesgue  measure,
and denote by the same letter  their density functions.
If   $p$ and $q$ are density functions on $\R^m$ then
$$
D(p\Vert q)
=
\int_{\R^m} p(u)\log \frac {p(u) } {q(u)}du\qute{if}\quad p<<q,
\tag 1.2
$$
and $\infty$ otherwise.
If $Z$ and $U$ are \rv s with values in $\R^m$ and
distributed according to $p=\Lw(Z)$ resp. $q=\Lw(U)$,
then we shall also use the notation
$D(Z\Vert U)$ for the  relative entropy $D(p\Vert q)$.
\enddefinition
\bsk

\definition{Definition}
For measures $p$ and $q$ on $\R^m$,
the Fisher information of  $p$ \wre\   $q$ is defined as
$$
I(p\Vert q )
=
\int_{\R^m}
\biggl|\nabla\log \frac {p(u)} {q(u)}\biggr|^2    p(du),
\tag 1.3
$$
if
$\log( {p(u)}/{q(u)})$ is smooth.
\enddefinition
\bsk

\definition{Definition}
The distribution
$q$ on $\R^m$ satisfies a logarithmic Sobolev \ineq\
with constant $\rho$ if
$$
D(p||q)\leq \frac 1{2\rho}\cdot I(p\Vert q)
$$
for all density functions $p$ on $\R^m$ with
$\log ({p(u)}/{q(u)})$  smooth.
\enddefinition
\bsk

A \LSI\   for a \pms\  $q$
is equivalent to the
hypercontractivity of the diffusion semigroup
 associated with  $q$.
 The prototype is  Gross'
\LSI\  for  Gaussian measure which is associated to the
Ornstein-Uhlenbeck semigroup [1], [2]. Another use of  logarithmic
Sobolev inequalities is to derive transportation cost
inequalities (a tool to prove measure concentration), c.f.
F. Otto, C. Villani [3]. The   logarithmic Sobolev \ineq\ for the stationary
distribution of a spin system
is equivalent to the property  called ``exponential decay
of correlation'';  for this concept we refer to
Bodineau and Helffer [4] and Helffer [5].
\bsk

In Euclidean spaces of dimension greater than $1$,
 no simple characterization is available for the measures $q$
satisfying a \LSI\ with some positive constant.
A well-known sufficient condition was given by Bakry and Emery [6]:
A density function $q(x)=\exp(-V(x))$ on $\R^m$
 satisfies a \LSI\ provided $V$ is
  uniformly  strictly convex.
Another useful result is Holley and Stroock's perturbation lemma [7]
which asserts that if   $q$ and $\tilde q$ are  density functions
on $\R^m$, such that the ratio $\tilde q(x)/ q(x)$ is bounded
both from above and below,
then $q$ and $\tilde q$ either both satisfy a \LSI, or neither of them does.

\bsk

For  measures on Euclidean spaces  with non-compact support,
it has been  a challenging task to derive \LSIs\
from  \LSIs\ for the local specifications.
(The  local specifications of the measure
$q=\Lw(X_1,\dots,X_m)$ on $\R^m$ are the conditional
densities $Q_i(\cdot|x_1,\dots, x_{i-1}, x_{i+1},\dots, x_m)
=\Lw(X_i|X_1=x_1,\dots, X_{i-1}= x_{i-1},X_{i+1}= x_{i+1},\dots,X_m= x_m)$.)
Let $q$ be a density function on a Euclidean space, and assume that the
local specifications of $q$ satisfy \LSIs\ with constants $\rho_i$.
It has been clear for a long time that
a reasonable approach to prove a \LSI\ for $q$
is to assume that
the mixed partial derivatives of $V(x)=-\log q(x)$
are not too large
relative to the numbers $\rho_i$.
This line was followed by
B. Zegarlinski [8] and, following in his footsteps,
G. Royer [9], Th\'eor\`{e}m 5.2.1).
Their results were improved by F. Otto and M. Reznikoff [10].
The present paper  follows this line, too.
The conditions of Otto and  Reznikoff's main theorem
helped to find the proper conditions
for the results in the present paper,
however,
our approach is entirely different from their's.
We shall discuss Otto and  Reznikoff's theorem at the end of Section 2.

\beginsection 2. {Statement of the results}

Let $\R^N$ denote the $N$-dimensional Euclidean  space
equipped with the Euclidean distance and the Borel
$\si$-algebra.
\bsk

Let us fix a density function
$$
q(x)=exp(-V(x)),\quad x\in \R^N.
$$
\bsk

We shall use the following
\bsk

{\bf Notation:}
\bigskip

$\bullet\quad q: \quad \text{ a fixed  density function   on}\quad  \R^N;$
\bsk

$\bullet\quad
X=(X_1,X_2,\dots,X_N): \quad\text{random \sqe\ in}\quad \R^N,\quad \Lw( X)=q;
$
\bsk

$\bullet\quad p:
\quad \text{another density function on}\quad \R^N; $
\bsk

$\bullet\quad
Y=(Y_1,Y_2,\dots,Y_N): \quad\text{random \sqe\ in}\quad \R^N,\quad \Lw( Y)=p;$
\bsk

$\bullet\quad(I_k, k=1,2,\dots,n)$: a partition of $[1,N]$, $|I_k|=n_k$;
\bsk

$\bullet\quad $ for  $x\in \R^N$,
$x^{(k)}\triangleq\bigl\{x_i: i\in I_k\bigr\},
\quad
\bar x^{(k)}\triangleq \bigl\{x_i: i\notin I_k\bigr\};$
\bsk

$\bullet\quad$
$X^{(k)}$ and $\bar X^{(k)}$: the corresponding
segments of $X$;
\bsk

$\bullet\quad$ $Y^{(k)}$ and $\bar Y^{(k)}$: the corresponding
segments of $Y$;
\bsk

$\bullet \quad
\bar q^{(k)}
\triangleq
\Lw\bigl(\bar X^{(k)}\bigr),\quad
  Q^{(k)}(\cdot|\bar x^{(k)})
\triangleq
\Lw\bigl(X^{(k)}|\bar X^{(k)}=\bar x^{(k)}\bigr);
$
\bsk

$
\bullet \quad
\bar p^{(k)}
\triangleq
\Lw\bigl(\bar Y^{(k)}\bigr),
 \qquad p^{(k)}(\cdot|\bar y^{(k)})
\triangleq
\Lw\bigl(Y^{(k)}|\bar Y^{(k)}=\bar y^{(k)}\bigr).
$

\bsk
\bsk
We consider $\R^N$ as the product of Euclidean spaces $\R^{(k)}$ of
dimension $n_k$.
\bsk

\definition{Definition}
The \cd s  $ Q^{(k)}(\cdot|\bar x^{(k)})$ and  $ p^{(k)}(\cdot|\bar x^{(k)})$ are
called the local specifications of $q$ resp. $p$.
\enddefinition
\bsk

To formulate the main results of this paper, we also need the concept
of (average) conditional relative entropy, together with some more notation:

\bsk

\definition{Definition}
If we are given a \pms\ $\pi=\Lw(S)$ on $\R^\ell$ ($\ell\ge 1$ integer), and
\cd s $\mu(\cdot|s)=\Lw(Z|S=s)$, $\nu(\cdot|s)=\Lw(U|S=s)$  on $\R^m$
then consider the average relative entropy
$$
\Bbb E_{\pi}D\bigl(\mu(\cdot|S)\Vert \nu(\cdot|S)\bigr)
=
\int_{\R^\ell}D\bigl(\mu(\cdot|s)\Vert \nu(\cdot|s)\bigr)\pi(ds).
$$
For $\Bbb E_{\pi}D(\mu(\cdot|S)\Vert \nu(\cdot|S))$
we shall use either of the notations
$$
 D\bigl(\mu(\cdot|S)\Vert \nu(\cdot|S)\bigr),
\quad D\bigl(\mu(\cdot|S)\Vert U|S\bigr),
\quad D\bigl(Z|S)\Vert \nu(\cdot|S)\bigr),
\quad D\bigl(Z|S)\Vert U|S)\bigr).
$$
\enddefinition
\bsk

For a fixed measure $q$ on $\R^N$, we want
to derive an inequality of the form
$$
D(p||q)
\leq
\frac 1 \rho\cdot
\sumk \rho_k\cdot D\bigl(p^{(k)}(\cdot|\bar Y^{(k)})||Q^{(k)}(\cdot|\bar Y^{(k)})\bigr)
\qute{for all}\quad p\qute{on}\quad \R^N,
\tag 2.1
$$
for some positive constants $\rho_k, 1\leq k\leq n$, and $\rho$.
I.e.,
we want to bound
$D(p\Vert q)$ by a weighted sum of the ``single phase''conditional
entropies
 $D\bigl(p^{(k)}(\cdot|\bar Y^{(k)})||Q^{(k)}(\cdot|\bar Y^{(k)})\bigr)$.
A bound of type (2.1)  holds only for a restricted class of \pms s $q$, and we want
a sufficient condition for (2.1).
Since relative entropy measures in a way
how different \pms s are,
\ineq\ (2.1) allows us to
 conclude to closeness of $p$ and $q$
 from the closeness of their local specifications.
Moreover,  an \ineq\ of type (2.1) ensures that
upper bounds  for the ``single phase''
relative entropies $D(p^{(k)}(\cdot|\bar  y^{(k)})||Q^{(k)}(\cdot|\bar y^{(k)}))$
that hold uniformly in $\bar y^{(k)}$, yield a bound for $D(p\Vert q)$.
This is a way to get \LSIs\ for measures on product spaces.
\bsk

To get \ineq\ (2.1), we make   three assumptions explained below.
Recall that $(I_k, k=1,2,\dots,n)$ is a partition of $[1,N]$.
\bsk

\definition{Assumption 1}
Assume that
$Q^{(k)}(\cdot|\bar x^{(k)})$ satisfies a \LSI\
 with constant $ \rho_k$
for all $x\in \R^N$ and  $k\in[1,n]$.
\enddefinition
\bsk

Consider the Hessian of $V(x)=-\log q(x)$, i.e., the matrix
$\bigl(V_{i,j}(x)\bigr)_{i,j\in [1,N]}$, where we denote by
$V_{i,j}(x)$  the second partial derivatives of $V(x)$.
\bsk

\definition{Assumption 2}
Assume that, for each $k\in[1,n]$,
the matrix $(V_{i,j}(x))_{i,j\in I_k}$
is bounded  from below
by some (possibly negative) constant times the identity.
\enddefinition
\bsk

To formulate Assumption 3, we introduce the following
\bsk

\definition{Notation}
Under Assumption 1, and for
\sqe s $x,\xi\in \R^N$ fixed, we denote
by  $A(x,\xi)$ the matrix   with elements
$$
\align
&A_{i,j}(x,\xi)
=
\frac {V_{i,j}(\bar x^{(\ell)},\xi^{(\ell)})}
{\sqrt{\rho_k}\cdot \sqrt{\rho_\ell}}
\qute{for}\quad i\in I_k, j\in I_\ell,k\neq \ell,\\
&A_{i,j}(x,\xi)=0 \qute{if}\quad i\qute{and}\quad j \qute{belong to the same set}\quad I_k.
\endalign
$$
Moreover, for  \sqe s $x,\xi\in \R^N$ and $0<\rho<\min \rho_k$,
we denote
by  $A^\rho(x,\xi)$ the matrix   with elements
$$
\align
&A^\rho_{i,j}(x,\xi)
=\frac {V_{i,j}(\bar x^{(\ell)},\xi^{(\ell)})}
{\sqrt{\rho_k-\rho}\cdot \sqrt{\rho_\ell-\rho}}
\qute{for}\quad i\in I_k, j\in I_\ell,k\neq \ell,\\
&A^\rho_{i,j}(x,\xi)=0 \qute{if}\quad i\qute{and}\quad j \qute{belong to the same set}\quad I_k.
\endalign
$$
(Thus $A(x,\xi)=A^0(x,\xi)$.)
\enddefinition
\bsk

\remark{Remark}
Unless the matrix  $A^\rho(x,\xi)$ is constant in $x$,
it is not symmetric, since
in the definition of $A^\rho_{i,j}(x,\xi)$ ($i\in I_k,j\in I_\ell$),
we use $\xi^{(\ell)}$, and not $\xi^{(k)}$.
\endremark

\definition{Assumption 3}
We assume that
$$
\sup_{x,\xi}\bigl\Vert A(x,\xi) \bigr\Vert\triangleq 1-\de<1,
\tag 2.2
$$
and that
$\rho$ is such that
$$
\sup_{x,\xi}\bigl\Vert A^\rho(x,\xi) \bigr\Vert\leq 1.
\tag 2.3
$$
\enddefinition
\bsk

Conditions (2.2) and (2.3) shall be used in the following form:
For all $x,\xi,u,v\in \R^N$,
$$
\align
&\Biggl|
\sum_{k,\ell\in[1,n],k\neq\ell}\quad
\sum_{i\in I_k,j\in I_\ell}
u_i\cdot V_{i,j} \bigl(x,\xi\bigr)\cdot v_j
\Biggr|\\
&\leq
(1-\de)\cdot \sqrt{\sumk \rho_k \cdot |u^{(k)}|^2}
\cdot
\sqrt{\sum_{\ell=1}^n \rho_\ell \cdot |v^{(\ell)}|^2}
\tag 2.4
\endalign
$$
and
$$
\align
&\Biggl|
\sum_{k,\ell\in[1,n],k\neq\ell}\quad
\sum_{i\in I_k,j\in I_\ell}
u_i\cdot V_{i,j} \bigl(x,\xi\bigr)\cdot v_j
\Biggr|\\
&\leq
\sqrt{\sumk (\rho_k-\rho) \cdot |u^{(k)}|^2}
\cdot
\sqrt{\sum_{\ell=1}^n (\rho_\ell-\rho) \cdot |v^{(\ell)}|^2},
\tag 2.5
\endalign
$$
respectively.

\proclaim{Theorem 1}
If Assumptions 1-3  hold  then
$$
D(p||q)
\leq
\frac 1 \rho\cdot
\sum_{k=1}^n \rho_k\cdot D\biggl(p^{(k)}\bigl(\cdot|\bar Y^{(k)}\bigr)
||Q^{(k)}(\cdot|\bar Y^{(k)})\biggr)
\tag 2.6
$$
for any \pms\ $p$ on $\R^N$.
\endproclaim
\bsk

\proclaim{Theorem 2}
Under Assumptions  1-3,  $q$ satisfies a \LSI\  with constant $\rho$.
\endproclaim
\bsk

Theorem 2 follows from Theorem 1, using Assumption 1 and
the fact that
by the definition of the  operator $\nabla$
$$
I(p\Vert q)
=
\sum_{k=1}^n
\Bbb E
I\bigl(p^{(k)}(\cdot|\bar Y^{(k)}) \bV Q^{(k)}(\cdot|\bar Y^{(k)})\bigr).
$$
\bsk

The statement of Theorem 2 was proved by F. Otto and M. Reznikoff [10],
under a condition similar to, but stronger than, Assumption 3.
 We discuss Otto and  Reznikoff's  theorem  at the end of this section.
 \bsk

Next we formulate a \LSI\ for a discrete time Markov process
governed by the local specifications $Q^{(k)}(\cdot|\bar y^{(k)})$.
\bsk

\definition{Definition of weighted Gibbs sampler}
\newline
Given  a partition $(I_k, k=1,2,\dots,n)$ of $[1,N]$, and local specifiations
$Q^{(k)}(\cdot|\bar y^{(k)}\bigr)$,
the weighted Gibbs sampler $\Gamma$ with weights $(\pi^{(k)}, k=1,2,\dots,n)$
is the Markov operator on the \pms s $p$ (on $\R^N$)
 defined by
$$
\Gamma =\sumk \pi_k \Gamma_k, \qquad
\Gamma_k(z|y) =\de\bigl(\bar y^{(k)}, \bar z^{(k)}\bigr)\cdot
Q^{(k)}\bigl(z^{(k)}|\bar y^{(k)}\bigr).
$$
(Here $\de$ denotes Kronecker's $\de$.)
\enddefinition
\bsk

\proclaim{Corollary to Theorem 1}
\newline
If Assumptions 1-3 hold  then
for the weighted Gibbs sampler $\Gamma$ with weights
$$
\bigl(\rho_k/ R, k=1,2,\dots,n\bigr),\qquad
R=\sum_k \rho_k,
$$
 we have
$$
 D(p||q)
\leq
\frac R \rho \cdot \biggl(D(p||q)-D(p\Gamma||q)\biggr).
\tag 2.7
$$
Thus
$$
D(p\Gamma^m||q)
\leq
\biggl(1-\frac  \rho R \biggr)^m\cdot D(p||q).
$$
\endproclaim
\bsk

(2.7) follows from Theorem 1 by the
\ineq\
$$
D(p\Gamma||q)\leq \frac 1 R\sumk \rho_k D(p\Gamma_k||q)
$$
(a consequence of the convexity of relative entropy)
and the identity
$$
D(p||q)- D(p\Gamma_k||q)
=
D\biggl(p^{(k)}\bigl(\cdot|\bar Y^{(k)}\bigr)
       ||Q^{(k)}(\cdot|\bar Y^{(k)})\biggr).
$$
(2.7) can be considered as a \LSI\ for the Gibbs sampler $\Gamma$.
Indeed, for the Markov process defined by $\Gamma$,
it bounds relative entropy (from  the stationary distribution)
by the decrease of relative entropy along the
Markov process.
\bsk

Next we formulate a transportation-cost \ineq\ that follows
from Theorem 2, using the Otto-Villani theorem (Theorem  1 in [3]). We need the following
definitions:
\bsk

\definition{Definition}
The quadratic Wasserstein distance between the
probability measures $r$ and  $s$ on  $\R^m$ is defined as
$$
W(r,s)  = \inf_\pi [E_\pi |\xi-\eta|^2]^{1/2},
$$
where $\xi$ and $\eta$ are \rv s with laws $r$ resp. $s$,
$|\xi-\eta|$ denotes Euclidean distance,
and infimum
is taken over all distributions $\pi=\Lw(\xi,\eta)$  with marginals $r$
and $s$.
\enddefinition
\bsk

\definition{Definition}
A \pms\ $s$ on $\R^m$ satisfies a transportation-cost \ineq\ with constant $\rho$
if
$$
W^2(r,s)\leq \frac 2\rho \cdot D(r\Vert s)
$$
for all \pms s $r$ on $\R^m$.
\enddefinition
\bsk

Transportation-cost \ineqs\ are useful in proving measure concentration \ineqs.
A transportation-cost \ineq\ for the case when $q$ is Gaussian,
was proved by Talagrand [11].
Otto and Villani  generalized   Talagrand's  \ineq\
 as follows:
\bsk

\proclaim{ Otto and Villani's theorem for Euclidean spaces} [3],[12]
\newline\
If a density function  on $\Bbb R^m$ satisfies a
logarithmic Sobolev \ineq\  then
it satisfies a transportation-cost \ineq\ with the same constant.
\endproclaim
\bsk

By Otto and Villani's theorem, Theorem 2 implies the following
\bsk

\proclaim{Theorem 3}
\newline
If Assumptions 1-3 hold  then $q$ satisfies
a transportation-cost \ineq\ with constant $\rho$.
\endproclaim
\bsk

In [13],  corrected in [14], the statement of Theorem 3, for equal $\rho_k$'s,
 was proved modulo an absolute constant factor.
\bsk

Now we compare Theorem 2  with the result of [10].
\bsk

In [10] the statement of Theorem 2 is proved
under the following condition in place of (2.3):
\bsk

For $k,\ell\in[1,n]$, $k\neq\ell$, and    $x\in \R^N$, consider the following minors of
the Hessian of $V(x)$:
$$
K_{k,\ell}(x)=\biggl(V_{i,j}(x)\biggr)_{i\in I_k,j\in I_\ell},
$$
and set
$$
\kappa_{k,\ell}
=
\sup_x
\bigl\Vert \bigl(K_{k,\ell}(x)\bigr) \bigr\Vert.
$$
Then  consider the $n\times n$ matrix
$$
K=\bigl(\kappa_{k,\ell}\bigr)_{k,\ell\in [1,n], k\neq \ell}.
$$
($K$ has $0$'s in the main diagonal.)
 Otto and  Reznikoff use the assumption that
$$
K\leq \Lambda \bigl(\{\rho_k-\rho\}\bigr),
\tag 2.3'
$$
where $\Lambda \bigl(\{\rho_k-\rho\}\bigr)$ denotes the $n\times n$
diagonal matrix  with elements
$\rho_k-\rho$.
With the notation
$$
\align
&{\kappa'}_{k,\ell}
=
\frac {\kappa_{k,\ell}}{\sqrt{\rho_k-\rho}\cdot \sqrt{\rho_\ell-\rho}}\\
&{K'}^\rho=\bigl({\kappa}'_{k,\ell}\bigr)_{k,\ell\in [1,n], k\neq \ell},
\endalign
$$
(2.3') can be written in the form
$$
{K'}^\rho\leq Id,
$$
where $Id$ is the $n\times n$ identity matrix.
Since ${K'}^\rho $ is symmetric,
this means that the largest  eigenvalue of ${K'}^\rho$ is $\leq 1$.
The elements of ${K'}^\rho $ are non-negative, thus, by Perron's theorem,
 the largest   eigenvalue of ${K'}^\rho$ equals $\bV{K'}^\rho\bV$. I.e.,
in [10] it is actually assumed that
$$
\bV{K'}^\rho\bV \leq 1,
\tag 2.3"
$$
which is clearly stronger  than (2.3).
\bsk

\remark{Remark}
\newline
If $q$ is Gaussian then  the Hessian of $V(x)$ does not depend on $x$.
Otto and Reznikoff's result is tight for Gaussian distributions
with attractive interactions. (For $\R^{(k)}=\R$;
attractivity means that $V_{i,j}\leq 0$ for $i\neq j$.)
For $q$ Gaussian and  $\R^{(k)}=\R$,
Theorem 2 can be formulated as follows:
If $\bV A^0\bV<1$ then
$q$ satisfies a \LSI\ with constant $\rho$, where
$\rho$ is the largest number satisfying
$$
\bV A^\rho\bV=1.
\tag 2.8
$$
Thus Theorem 2 is tight for those  Gaussian distributions $q$ for which
$ \bV A^\rho\bV$ (for the $\rho$ defined by (2.8)) is given by
the absolute value of the smallest negative
eigenvalue  (and not the largest positive one).
\endremark
\bsk

\remark{Example}
\newline
Assumption 3 is practically impossible to check, except
when the mixed partial derivatives of $V(x)$ are constants.
Otherwise we probably cannot do better
than use Otto and Reznikoff's theorem.
However, if the mixed partial derivatives of $V(x)$ are all constants
then Theorem 2 may give a better result.
Indeed, let
$V(x)=-\log q(x)$ be of the  form
$$
V(x)=\sum_{k=1}^n \phi_k(x)+\sum_{k,\ell\in[1,n],
k\neq\ell} a_{k,\ell}\cdot x_k\cdot x_{\ell},
$$
where for each $k$ and fixed $\bar x_k$,
the single phase density $C_k(\bar x_k)\cdot \exp(-\phi_k(x_k,\bar x_k))$,
as a function of $x_k$,
satisfies a \LSI\ with a common constant $\rho$.
Theorem 2
guaranties a positive logarithmic Sobolev constant if
the matrix with elements
$$
 a_{k,\ell}\qute{outside the main diagonal, and}\quad 0 \qute{otherwise}
$$
 has norm $<\rho$.  On the other hand,
Otto and Reznikoff's theorem guaranties
a positive logarithmic Sobolev constant
if the matrix with elements
$$
| a_{k,\ell}|\qute{outside the main diagonal, and}\quad 0 \qute{otherwise}
$$
 has norm $<\rho$.
 To see a concrete example when the first condition holds,
but the second does not, consider the infinite dimensional Toeplitz matrix
$B=(b_{k,\ell})$ defined by
$$
 b_{k,k+1}= 1, \quad  b_{k,k+2}=-1,
\quad  b_{k,\ell}=0\qute{for}\quad \ell\ge k,\quad \ell\notin \{k+1,k+2\},
$$
and $b_{k,\ell}=b_{\ell,k}$.
From the theory of Toeplitz matrices (c.f. [15])
we know that
$$
\bV B\bV=2\cdot\max\bigl| \cos x-\cos(2x)\bigr|=\frac 9 4,
$$
while for the matrix $abs(B)$ consisting of the absolute values of
$b_{k,\ell}$, we get
$$
\bV abs(B)\bV =2\cdot\max\bigl| \cos x+\cos(2x)\bigr|=4.
$$
Denote by $B_m$ and  $abs(B_m)$  the matrices consisting of the first
 $m$ rows and columns of $B$ resp. $abs(B)$; clearly
$|| B_m||\leq  9 /4$ and $\lim_{m\to\infty}||abs(B_m)||=4$.
Therefore, if we take $A=B_m$, and
if the functions $\phi_k(x)$ in the definition of $V(x)$
are such that  the single phase densities
$C_k(\bar x_k)\cdot \exp(-\phi_k(x_k,\bar x_k))$
satisfy a \LSI\ with a common constant
$>\frac 9 4$ then Theorem 2
guaranties a positive logarithmic Sobolev constant for $q=\exp(-V)$.
However, we cannot get this
from Otto and Reznikoff's theorem.
\endremark
\bsk

\beginsection 3. {Proof of Theorem 1.}

Our approach to prove Theorem 1 is based on the
interpolation between the \pms s $p$ and $q$
realized by the
solution
of the Fokker-Planck equation
$$
\partial_t p_t(y)=\triangle  p_t(y)+\nabla\cdot\bigl(p_t(y)\cdot\nabla V(y)\bigr),
  \qquad p_0(y)=p(y).
\tag 3.1
$$
With the notation
$$
h=p/{q} \qute{and}\quad h_t=p_t/q,
$$
the Fokker-Planck equation (3.1) can be rewritten  as follows:
$$
\partial_t h_t
=
Lh_t\triangleq
\triangle h_t-\nabla h_t\cdot \nabla V,\qquad \quad h_0(y)= h(y).
\tag 3.2
$$
We have
$$
\Bbb E_p \log h=D(p||q)\qute{and}\quad \Bbb E_{p_t}\log h_t=D(p_t|| q).
$$
\bsk

Our argument heavily draws on the ideas developed
in the paper by F. Otto and C. Villani [3].
To be able to use the tools
of [3], we need the limit relation
$$
\lim_{t\to\infty}
D\bigl(p_t\bV q\bigr)=0.
\tag 3.3
$$
To this end we
 prove a \LSI\ for $q$ with a much smaller constant
than  claimed in Theorem 2. (It is disturbing that this weak
preliminary result requires a
very lengthy proof.)
\bsk

\proclaim{Auxiliary Theorem}
If Assumptions 1-3 hold  then $q$ satisfies a \LSI\  with a constant
$C=C(R,\rho_{min},\de)$, where $R=\sumk\rho_k$, $\rho_{min}=\min_k \rho_k$ and
$\delta=1-\sup_{x,\xi}|| A(x,\xi)||$.
\endproclaim
\bsk

For the proof of Theorem 1 we also need the following simple
lemma ( c. f. (32) in [3]).
\bsk

\proclaim{Approximation Lemma}
In the proof of Theorem 1 we can restrict ourselves to the case when
$V(x)=-\log q(x)\in \Cal C^\infty$, and
$h(x)=p(x)/q(x)$ is of the form
$$
\align
&h(x)=(1-\e)\cdot g(x)+\e,\qute{where}\\
& g\in\Cal C^\infty\qute{is a  compactly supported density function (\wre\ } q), \qute
{and}\quad \e >0.
\endalign
$$
\endproclaim
\bsk

The proofs of the Auxiliary Theorem and the Approximation Lemma
are postponed to Section 4,
 although they are used
in the proof of Theorem 1 in this section.
\bsk

We need  some more
\bsk

\definition{Notation}
\newline
Let
$$
Y_t=(Y_{t,1},Y_{t,2},\dots,Y_{t,N})
$$
denote a random \sqe\ with
$\Lw(Y_t)=p_t$, where $p_t$ is the solution of the Fokker-Planck
equation (3.1).
In accordance with  the notation at the beginning of Section 2, we write
$$
 Y_t^{(k)}
=
\bigl\{Y_{t,i}: i\in I_k\bigr\},
\quad
 \bar  Y_t^{(k)} =\bigl\{Y_{t,i}: i\notin I_k \bigr\}.
$$
Further, we set
$$
\bar p_t^{(k)}
=
\Lw\bigl(\bar Y_t^{(k)} \bigr),
\quad
p_t^{(k)} \bigl(\cdot|\bar y_t^{(k)}  \bigr)
=
\Lw\bigl(  Y_t^{(k)}|\bar  Y_t^{(k)}=\bar      y_t^{(k)}\bigr).
$$
\enddefinition
\bsk

By the Approximation Lemma we may assume that $V\in\Cal C^\infty$.  Then
 the domain of the operator $L$ in (3.2) can be defined
 so as to contain the class  $\Cal D_0$
of those functions $h$ in $\Cal C^\infty$ that are bounded, and
 whose partial derivatives of any order,
 multiplied by the partial derivatives of $V$ of any order,
are bounded.
The class  $\Cal D_0$ is dense in $\Bbb L_2(q)$ and  stable under $L$.
\bsk

Again by the Approximation Lemma
we can assume that the function $h_0=h$ in (3.2) belongs to $\Cal D_0$.
As explained in [3], this implies that
$h_t$ is uniformly bounded from below and from above, and that, for $t$ fixed,
$|\nabla h_t|^2$ is bounded. (Here we use the fact that, by Assumptions 2 and 3,
the Hessian of $V(x)$ is bounded from below
by a (possibly negative) constant times the identity.--
In [3]   assumption (32) of that paper is used which is implied by
the  assumption $h_0=h\in\Cal D_0$.)
\bsk

Consequently, as explained in [3], under condition $h_0=h\in\Cal D_0$,
the Fokker-Planck equation (3.2)
defines a semigroup of diffeomorphisms
$$
\Phi_t: \R^N\mapsto \R^N,\quad 0\leq t<\infty,
\tag 3.4
$$
satisfying
$$
\partial_t \Phi_t(y)=-\nabla \log h_t(\Phi_t(y)),
\tag 3.5
$$
and
$$
p_t\Phi_s=p_{t+s}.
\tag 3.6
$$
(3.6) means that
$p_{t+s}$ is the image of $p_t$ under the map $\Phi_s$.
Since $\Lw(Y_t)=p_t$,
we can think of the random \sqe s $Y_t$ as
functions of $Y=Y_0$:
$$
Y_t=\Phi_t(Y)=\Phi_t(Y_0).
$$
\bsk


Let us introduce the function
$$
\chi(y)=\sumk \rho_k\bigl[\log h(y)-\log \bar h^{(k)}(\bar y^{(k)})\bigr],\qquad y\in \R^N,
$$
where
$$
\bar h^{(k)}\bigl(\bar  y^{(k)}\bigr)
=
\frac
{\bar p^{(k)}\bigl(\bar y^{(k)}\bigr)}
{\bar q^{(k)}\bigl(\bar y^{(k)}\bigr)}
=
\int_{\R^{(k)}}h(y)Q^{(k)}
\bigl(dy^{(k)})|\bar y^{(k)}\bigr).
$$
(The integration domain is $\R^{n_k}$; the superscript $(k)$ indicates
that integration is \wre\ the  variable $y^{(k)}$.)
We  have
$$
\Bbb E_p \chi
=
\sumk \rho_k
\cdot
D\biggl(p^{(k)}\bigl(\cdot|\bar Y^{(k)}\bigr)
\bV Q^{(k)}\bigl(\cdot|\bar Y^{(k)}\bigl)\biggr).
$$
Thus the statement of Theorem 1 is equivalent to
$$
\rho\cdot\Bbb E_p\log h
\leq
\Bbb E_p \chi.
$$
\bsk

It is well known (and a proof can be found in [3]) that
$$
\frac{\partial}{\partial t} D(p_t|| q)
=
-I(p_t|| q)
=
-\Bbb E_{p_t}\bigl|\nabla \log h_t\bigr|^2.
$$
Thus, by (3.3),
$$
D(p|| q)=D(p|| q)-\lim_{t\to\infty}D(p_t|| q)
=
\int_0^\infty
\Bbb E_{p_t}\bigl|\nabla \log h_t\bigr|^2 dt.
\tag 3.7
$$
\bsk

We introduce, analogously to the definition of $\chi$,
the functions
$$
\chi_t(y)=
\sumk\rho_k\bigl[\log h_t(y)-\log \bar h_t^{(k)}\bigl(\bar y^{(k)}\bigr)\bigr],
$$
where
$$
\bar h_t^{(k)}\bigl(\bar y^{(k)}\bigr)
=
\bar p_t^{(k)}\bigl(\bar y^{(k)}\bigr)
/{\bar q^{(k)}\bigl(\bar y^{(k)}\bigr)}.
$$
We have
$$
\Bbb E_{p_t} \chi_t
=
\sumk \rho_k
\cdot
D\biggl(p_t^{(k)}\bigl(\cdot|\bar Y_t^{(k)}\bigr)
\bV Q^{(k)}\bigl(\cdot|\bar Y_t^{(k)}\bigl)\biggr).
$$
In particular, $\Bbb E_{p_t} \chi_t\ge 0$.

\bsk

Using (3.7) and the fact that
$\Bbb E_{p_t} \chi_t\ge 0$, for the proof of
Theorem 1 it is enough to prove the following two propositions:
\bsk

\proclaim{Proposition 1}
$$
\Bbb E_p\chi
-
\lim_{t\to\infty}\Bbb E_{p_t} \chi_t
=
\int_0^\infty
\Bbb E_{p_t} \biggl\{\nabla\chi_t\cdot \nabla\log h_t\biggr\} dt.
\tag 3.8
$$
\endproclaim
\bsk

\proclaim{Proposition 2}
$$
\Bbb E_{p_t}\biggl\{\nabla\chi_t\cdot \nabla\log h_t\biggr\}
\ge
\rho\cdot
\Bbb E_{p_t} \bigl|\nabla\log h_t\bigr|^2.
$$
\endproclaim
\bsk

\demo{Proof of Proposition 1}
\newline
For all $y\in\R^N$ we have
$$
\chi(y)-\lim_{t\to\infty} \chi_t\bigl(\Phi_t(y)\bigr)
=
-\int_0^\infty
\frac {\partial} {\partial t}\bigl( \chi_t\bigl(\Phi_t(y)\bigr)\bigr)dt.
$$
Therefore, by Fubini's theorem,
$$
\Bbb E_p
\biggl\{\chi(Y)-\lim_{t\to\infty} \chi_t\bigl(\Phi_t(Y)\bigr)\biggr\}
=
-\int_0^\infty
\Bbb E_p
\biggl\{\frac {\partial} {\partial t} \bigl(\chi_t(\Phi_t(Y))\bigr)\biggr\}dt.
\tag 3.9
$$
\bsk

Denoting by dot derivation \wre\ $t$,  and using (3.5):
$$
\frac {\partial} {\partial t}\biggl( \chi_t\bigl(\Phi_t(y)\bigr)\biggr)
=
\dot\chi_t\bigl(\Phi_t(y)\bigr)
-
\nabla\chi_t\bigl(\Phi_t(y)\bigr)\cdot\nabla\log h_t\bigl(\Phi_t(y)\bigr).
\tag 3.10
$$
\bsk

Further,
$$
\align
&\partial_t\bigl(\chi_t(z)\bigr)
=
\sumk
\rho_k\cdot
\biggl[\partial_t\biggl(\log h_t(z)\biggr)
-
\partial_t
\biggl(\log\bar h_t^{(k)}\bigl(\bar z^{(k)}\bigr) \biggr)
\biggr],   \quad  z\in \R^N.
\tag 3.11
\endalign
$$
\bsk

To calculate $\partial_t(\log\bar h_t^{(k)}\bigl(\bar z^{(k)}\bigr))$,
we need the following
\bsk

\proclaim{Lemma}
\newline
The solution $h_t$ of the Fokker-Planck equation (3.2) satisfies
$$
||\partial_t h_t||_{\Bbb L_2(q)}\leq ||Lh_0||_{\Bbb L_2(q)}<\infty.
\tag 3.12
$$
\endproclaim
\bsk

\demo{Proof}
\newline
The operator $L$ is
defined on a dense subset $\Cal D_0$ of $\Bbb L_2(q)$.
Moreover, $L$ is symmetric and negative definite
on $\Cal D_0$. Indeed, by partial integration we have
$$
(Lf,g)_{L_2(q)}=\int_{\R^N}\bigl(\triangle f-\nabla V \cdot\nabla f\bigl)\cdot g dq
=
-\int_{\R^N}\nabla f\cdot\nabla g dq.
$$
It follows that for $\lambda>0$
$$
\bigl((\lambda I-L)f,f\bigr)_{\Bbb L_2(q)}\ge \lambda \cdot||f||^2_{\Bbb L_2(q)},
$$
i.e.,
$$
\bV\bigl(\lambda I-L\bigr)^{-1}\bV \leq \frac 1 \lambda.
$$
Thus by the Hille-Yosida theorem  (c.f. [16]),
there exists a contraction semigroup $(P_t:t\ge0)$ on $\Bbb L_2(q)$
whose generator is $L$:
$$
\partial_t P_th_0=L P_th_0\qute{for}\quad  h_0\in  \Cal D_0,\qute{and}\quad ||P_t||\leq 1.
$$
For $h_0\in \Cal D_0$,
the solution of (3.2) can be written as
$h_t=P_th_0$,
and since $P_tL=LP_t$, we have
$$
\partial_t h_t
=
\partial_tP_th_0=LP_th_0=P_tLh_0,
$$
which implies (3.12). $\qquad\qquad\qquad\qquad\qquad\qquad\qquad\qquad\qquad\qquad\qquad\qquad\quad\qed$

\bsk

\enddemo

By the above Lemma, $\partial_t h_t \in \Bbb L_1(q)$,
 so we can differentiate under the integral sign
in the next formula:
$$
\align
&\partial_t
\biggl(\log \bar h_t^{(k)}(\bar z^{(k)}) \biggr)
=
\partial_t
\int_{\R^{(k)}}  h_t(z)
Q^{(k)}\bigl(dz^{(k)}|\bar z^{(k)}\bigr)\\
&=
\frac
{\int_{\R^{(k)}}
\partial_t \biggl( h_t(z)    \biggr)
Q^{(k)}\bigl(dz^{(k)}|\bar z^{(k)}\bigr)}
{\bar h_t^{(k)} \bigl(\bar z^{(k)}\bigr)}\\
&=
\int_{\R^{(k)}}
\partial_t
\log h_t(z) \cdot
\frac  {h_t(z)}
{\bar h_t^{(k)}
  \bigl(\bar  z^{(k)}\bigr)}
Q^{(k)}\bigl(dz^{(k)}|\bar z^{(k)}\bigr).
\tag 3.13
\endalign
$$
\bsk

By the definition of the function $h_t$,
$$
\frac  {h_t(z)}{\bar h_t^{(k)}   \bigl(\bar  z^{(k)})}
Q^{(k)}\bigl(dz^{(k)}|\bar z^{(k)}\bigr)
=
p_t^{(k)}\bigl(dz^{(k)}|\bar z^{(k)}\bigr).
$$
Thus (3.13) implies
$$
\align&
\partial_t
\biggl(\log\bar h_t^{(k)}\bigl(\bar z^{(k)}\bigr)\biggr)
\\&
=
\int_{\R^{(k)}}
\partial_t
\log h_t(z)p_t^{(k)}\bigl(dz^{(k)}|\bar z^{(k)}\bigr)
=
\Bbb E_{p_t}\bigl\{\partial_t\log h_t|\bar z^{(k)} \bigr\},
\tag 3.14
\endalign
$$
where $\bar z^{(k)}$ in the condition of the  expectation
is a shorthand for $\bar Y_t^{(k)}=\bar z^{(k)}$.
Substituting (3.14) into (3.11) we get
$$
\partial_t
\bigl(\chi_t(z)\bigr)
=
\sumk
\rho_k\cdot
\biggl[\partial_t\log h_t(z)
-
\Bbb E_{p_t}\bigl\{\partial_t\log h_t|\bar z^{(k)}   \bigr\}
\biggr].
$$
It follows that
$\Bbb E_{p_t} \dot \chi_t=0$
which, together with (3.10), yields
$$
\Bbb E_p
\frac {\partial} {\partial t}\biggl( \chi_t\bigl(\Phi_t(Y)\bigr)\biggr)
=
-\Bbb E_{p_t} \biggl\{\nabla\chi_t\cdot \nabla\log h_t\biggr\}.
$$
Substituting this into (3.9) we get (3.8).
$\qquad\qquad\qquad\qquad\qquad\qquad\qquad\qquad\qquad\qed$

\enddemo
\bsk

\demo{Proof of Proposition 2}
\newline
We prove Proposition 2 for $t=0$; for $t>0$ the proof is the same.
For a function  $g:\R^N\mapsto \R$ set
$$
\nabla^{(k)}g(x)=\bigl(\partial_i g(x): i\in I_k\bigr).
$$

\bsk

We need the following
\bsk

\proclaim{Proposition 3}
\newline
For $k, \ell\in[1,n]$,   $ k\neq \ell$, we have
$$
\align
& \nabla^{(k)}  \log \bar h^{(\ell)}\bigl(\bar y^{(\ell)}\bigr)
=
\Bbb E_p\bigl\{\nabla^{(k)} \log h|\bar y^{(\ell)}\bigr\}\\
&-
\int_{\R^{(\ell)}\times \R^{(\ell)}}
\bigl[
\nabla^{(k)}V(\bar y^{(\ell)},\xi^{(\ell)})
-
\nabla^{(k)}V(\bar y^{(\ell)},\eta^{(\ell)})
\bigr]
\pi^{(\ell)}(d\xi^{(\ell)},d\eta^{(\ell)}|\bar y^{(\ell)}),
\tag 3.15
\endalign
$$
where $\pi^{(\ell)}(d\xi^{(\ell)},d\eta^{(\ell)}|\bar y^{(\ell)})$
 is an  arbitrary coupling of
the conditional measures\newline
$p^{(\ell)}\bigl(\cdot|\bar y^{(\ell)}\bigr)$
and
$Q^{(\ell)}\bigl(\cdot|\bar y^{(\ell)}\bigr)$.
  (I.e.,  $\pi^{(\ell)}(d\xi^{(\ell)},d\eta^{(\ell)}|\bar y^{(\ell)})$
is a  conditional
density on  $\R^{(\ell)}\times \R^{(\ell)}$           with  marginals
$p^{(\ell)}(\cdot|\bar y^{(\ell)})$
and
$Q^{(\ell)}(\cdot|\bar y^{(\ell)})$.)
\endproclaim
\bsk

\demo{Proof of Proposition 3}
\newline
Since $|\nabla h|$ is bounded (and $|\nabla h_t|$ is also bounded for $t$ fixed),
we have
$$
\nabla^{(k)}   \bar h^{(\ell)}(\bar y^{(\ell)})
=
\int_{\R^{(\ell)}}
\nabla^{(k)}
 \biggl(h\bigl(\bar y^{(\ell)}, \xi^{(\ell)}\bigr)\cdot
 Q^{(\ell)}\bigl(\xi^{(\ell)}|\bar y^{(\ell)}\bigr)\biggr) d\xi^{(\ell)}.
\tag 3.16
$$
Further,
$$
\align
&\nabla^{(k)}  Q^{(\ell)}\bigl(\xi^{(\ell)}|\bar y^{(\ell)}\bigr)
=
\nabla^{(k)}
\frac{\exp\bigl(-V\bigl( \bar y^{(\ell)}, \xi^{(\ell)}) \bigr) }
     {\int_{\R^{(\ell)}}
     \exp\bigl(-V(\bar y^{(\ell)},\eta^{(\ell)})\bigr) d \eta^{(\ell)}}\\
&=
-
\nabla^{(k)}
V\bigl(\bar y^{(\ell)}, \xi^{(\ell)}\bigr)
\cdot Q^{(\ell)}( \xi^{(\ell)}| \bar y^{(\ell)})\\
&\qquad\qquad\qquad+
Q^{(\ell)}( \xi^{(\ell)}| \bar y^{(\ell)})
 \cdot
 \int_{\R^{(\ell)}}  \nabla^{(k)}V( \bar y^{(\ell)},\eta^{(\ell)})
  Q^{(\ell)}(d\eta^{(\ell)}| \bar y^{(\ell)})\\
&=
Q^{(\ell)}( \xi^{(\ell)}| \bar y^{(\ell)})
 \cdot
  \int_{\R^{(\ell)}}
  \biggl[ \nabla^{(k)}V( \bar y^{(\ell)},\eta^{(\ell)})
  -
  \nabla^{(k)}V\bigl(\bar y^{(\ell)}, \xi^{(\ell)}\bigr)   \biggr]
  Q^{(\ell)}(d\eta^{(\ell)}| \bar y^{(\ell)}).
\endalign
$$
\bsk

It follows that
$$
\align
&\nabla^{(k)} \biggl(h\bigl(\bar y^{(\ell)}, \xi^{(\ell)}\bigr)
\cdot
Q^{(\ell)}\bigl( \xi^{(\ell)}| \bar y^{(\ell)}\bigr)\bigr)\biggr)\\
&=
Q^{(\ell)}\bigl( \xi^{(\ell)}| \bar y^{(\ell)}\bigr)\bigr)
\cdot
\Biggl[
 \nabla^{(k)}h\bigl(\bar y^{(\ell)}, \xi^{(\ell)}\bigr)\\
&+
h\bigl(\bar y^{(\ell)}, \xi^{(\ell)}\bigr)
\cdot
\int_{\R^{(\ell)}}
\biggl(
\nabla^{(k)}V( \bar y^{(\ell)},\eta^{(\ell)})
-
\nabla^{(k)}V \bigl(\bar y^{(\ell)}, \xi^{(\ell)}\bigr)  \biggr)
   Q^{(\ell)}\bigl( d\eta^{(\ell)}| \bar y^{(\ell)}\bigr)\Biggr].
\tag 3.17
\endalign
$$
Substituting  (3.17) into (3.16):
$$
\align
& \nabla^{(k)}   \bar h^{(\ell)}(\bar y^{(\ell)})\\
&=
\int_{\R^{(\ell)}}
\nabla^{(k)} h\bigl(\bar y^{(\ell)}, \xi^{(\ell)}\bigr)
   Q^{(\ell)}(d\xi^{(\ell)}| \bar y^{(\ell)})\\
&+
\int_{\R^{(\ell)}}
  h\bigl(\bar y^{(\ell)}, \xi^{(\ell)}\bigr)   Q^{(\ell)}(d\xi^{(\ell)}| \bar y^{(\ell)})
\cdot
\int_{\R^{(\ell)}}
\nabla^{(k)}V( \bar y^{(\ell)},\eta^{(\ell)}) Q^{(\ell)}(d\eta^{(\ell)}| \bar y^{(\ell)}) \\
&-
\int_{\R^{(\ell)}}
  h\bigl(\bar y^{(\ell)}, \xi^{(\ell)}\bigr)
  \cdot \nabla^{(k)}V \bigl(\bar y^{(\ell)}, \xi^{(\ell)}\bigr)
  Q^{(\ell)}\bigl( d\xi^{(\ell)}| \bar y^{(\ell)}\bigr )  \\
&=
\int_{\R^{(\ell)}}
\nabla^{(k)} h\bigl(\bar y^{(\ell)}, \xi^{(\ell)}\bigr)
Q^{(\ell)}(d\xi^{(\ell)}| \bar y^{(\ell)})
+
\bar h^{(\ell)} (\bar y^{(\ell)})
\cdot   \int_{\R^{(\ell)}}
\nabla^{(k)}V( \bar y^{(\ell)},\eta^{(\ell)}) Q^{(\ell)}(d\eta^{(\ell)}| \bar y^{(\ell)})\\
&-
\int_{\R^{(\ell)}}
 h\bigl(\bar y^{(\ell)}, \xi^{(\ell)}\bigr)
\cdot \nabla^{(k)}V \bigl(\bar y^{(\ell)}, \xi^{(\ell)}\bigr)
   Q^{(\ell)}(d\xi^{(\ell)}| \bar y^{(\ell)}).
\endalign
$$
Dividing both sides by $ \bar h^{(\ell)}(\bar y^{(\ell)} ) $:
$$
\align
& \nabla^{(k)}  \log \bar h^{(\ell)}(\bar y^{(\ell)})
=
\int_{\R^{(\ell)}}
\nabla^{(k)}\log h\bigl(\bar y^{(\ell)}, \xi^{(\ell)}\bigr)
\cdot \frac {h\bigl(\bar y^{(\ell)}, \xi^{(\ell)}\bigr)}
{\bar h^{(\ell)}\bigl(\bar y^{(\ell)}\bigr)}
Q^{(\ell)}(d\xi^{(\ell)}| \bar y^{(\ell)})\\
&+
\int_{\R^{(\ell)}}
\nabla^{(k)}V\bigl( \bar y^{(\ell)},\eta^{(\ell)}\bigr)
Q^{(\ell)}(d\eta^{(\ell)}| \bar y^{(\ell)})\\
&-
\int_{\R^{(\ell)}}
\nabla^{(k)}V\bigl(\bar y^{(\ell)}, \xi^{(\ell)}\bigr)
\cdot \frac {h\bigl(\bar y^{(\ell)}, \xi^{(\ell)}\bigr)}
{\bar h^{(\ell)}\bigl(\bar y^{(\ell)}\bigr)}
Q^{(\ell)}(d\xi^{(\ell)}| \bar y^{(\ell)}).
\endalign
$$
Since
$$
\frac {h\bigl(\bar y^{(\ell)}, \xi^{(\ell)}\bigr)}
{\bar h^{(\ell)}\bigl(\bar y^{(\ell)}\bigr)}
Q^{(\ell)}(d\xi^{(\ell)}| \bar y^{(\ell)})
=
p^{(\ell)}(d\xi^{(\ell)}| \bar y^{(\ell)}),
$$
and
$$
\int_{\R^{(\ell)}}
\nabla^{(k)}\log h\bigl(\bar y^{(\ell)}, \xi^{(\ell)}\bigr)
 p^{(\ell)}(d\xi^{(\ell)}| \bar y^{(\ell)})
=
\Bbb E_p\bigl\{\nabla^{(k)}\log h |\bar y^{(\ell)}\bigr\},
$$
(3.15) follows. $\qquad\qquad\qquad\qquad\qquad\qquad\qquad\qquad\qquad\qquad\qquad\qquad\qquad\qquad\qed$

\enddemo
\bsk
\bsk

Now we are ready to prove Proposition 2. By Proposition 3 we have
$$
\align
&\nabla^{(k)} \chi(y)
=
\sum_{\ell=1}^n
\rho_\ell\cdot
\biggl[ \nabla^{(k)}\log h(y)
-
\nabla^{(k)}\log \bar h^{(\ell)}(y^{(\ell)})\biggr]\\
&=
\sum_{\ell=1}^n
\rho_\ell\cdot
\biggl[ \nabla^{(k)}\log h(y)
-\Bbb E_p\bigl\{ \nabla^{(k)} \log \bar h|\bar y^{(\ell)}\bigr\}\biggr]\\
&+
\int_{\R^N\times \R^N}
\biggl(\sum_{\ell\neq k}
\rho_\ell\cdot
\biggl[\nabla^{(k)} V\bigl(\bar y^{(\ell)}, \xi^{(\ell)}\bigr)
- \nabla^{(k)} V\bigl(\bar y^{(\ell)}  ,\eta^{(\ell)}\bigr)\biggr]\biggr)
\Pi(d\xi,d\eta|y),
\endalign
$$
where $\Pi(d\xi,d\eta|y)$ denotes the conditional product measure
$\prod_{\ell=1}^n \pi^{(\ell)}\bigl(d\xi^{(\ell)},d\eta^{(\ell)}|\bar y^{(\ell)}\bigr)$.
\bsk

It follows that
$$
\align
&\Bbb E_p\bigl\{\nabla \chi\cdot \nabla \log h\bigr\}
=
\sum_{k=1}^n \rho_k\cdot \Bbb E_p |\nabla^{(k)} \log h|^2\\
&+
\sum_{k,\ell\in [1,n], \ell\neq k}
\rho_\ell\cdot
\biggl[\Bbb E_p | \nabla^{(k)}  \log h|^2
-\Bbb E_p\biggl\{\Bbb E_p\bigl\{ \nabla^{(k)} \log h|\bar y^{(\ell)}\bigr\}
\cdot  \nabla^{(k)} \log h\biggr\}\biggr]\\
&+
\Bbb E_{p,\Pi}
\biggl\{
\sum_{k,\ell\in [1,n], \ell\neq k}
\rho_\ell\cdot
\bigl[ \nabla^{(k)}V\bigl(\bar y^{(\ell)}, \xi^{(\ell)}\bigr)
 -\nabla^{(k)}V\bigl(\bar y^{(\ell)},\eta^{(\ell)}\bigr)\bigr]
\cdot  \nabla^{(k)}  \log h.
\biggr\},
\tag 3.18
\endalign
$$
Here $\Bbb E_{p,\Pi}$ denotes expectation \wre\ the \jd\ $\Lw(Y,\xi,\eta)$, defined by
$\Lw(Y)=p$ and
$\Lw(\xi,\eta|Y)=\Pi(d\xi,d\eta|y)$.
\bsk

For  $k\neq \ell$ we have
$$
\align
&\Bbb E_p
\biggl\{
\biggl[
\nabla^{(k)} \log h(y)-\Bbb E_p\bigl\{\nabla^{(k)} \log h|\bar y^{(\ell)}\bigr\}\biggr]\cdot
\nabla^{(k)} \log h(y)\biggr\}\\
&=
\Bbb E_p\biggl|\nabla^{(k)}\log h(y)\biggr|^2
-
\Bbb E_p\biggl\{\Bbb E_p^2\bigl\{\nabla^{(k)} \log h|\bar y^{(\ell)}\bigr\}\biggr\}
\ge 0.
\tag 3.19
\endalign
$$
\bsk

To estimate the last line in (3.18),
we introduce the notation
$$
U(y,\xi)
=
\sum_{k,\ell\in [1,n], \ell\neq k}
\rho_\ell\cdot
\nabla^{(k)}V\bigl(\bar y^{(\ell)},\xi^{(\ell)}\bigr)\cdot \nabla^{(k)} \log h(y),
\quad y,\xi\in \R^N.
$$
We have
$$
\sum_{k,\ell\in [1,n], \ell\neq k}
\rho_\ell\cdot
\biggl[\nabla^{(k)}V\bigl(\bar y^{(\ell)}, \xi^{(\ell)}\bigr)
- \nabla^{(k)} V\bigl(\bar y^{(\ell)},\eta^{(\ell)}\bigr)\biggr]
\cdot \nabla^{(k)} \log h(y)
=
U(y,\xi)-U(y,\eta).
$$
\bsk

To estimate $|U(y,\xi)-U(y,\eta)|$, we carry out the followong calculation:
\bsk

$$
\align
&\frac{\partial} {\partial \tau} U\bigl(y,\eta+\tau(\xi-\eta)\bigr)\\
&=
\sum_{k,\ell\in [1,n], \ell\neq k}\quad \sum_{i\in I_k,j\in I_\ell}
\rho_\ell\cdot(\xi_j-\eta_j)
\cdot V_{i,j}\bigl(\bar y^{(\ell)}, \eta^{(\ell)}+\tau(\xi^{(\ell)}-\eta^{(\ell)})\bigr)
\cdot \partial_i \log h(y).
\endalign
$$
\bsk

Hence, by Assumption 3 (c.f. (2.5)),
$$
\align
&\biggl|\frac{\partial} {\partial \tau} U\bigl( y,\eta+\tau(\xi-\eta)\bigr)\biggr|\\
&\leq
\sqrt{\sum_{\ell=1}^n\sum_{j\in I_\ell}(\rho_\ell-\rho)\cdot \rho_\ell^2\cdot (\xi_j-\eta_j)^2}
\cdot
\sqrt{\sum_{k=1}^n \sum_{i\in I_k}
(\rho_k-\rho)\cdot |\partial_i \log h(y)|^2}\\
&=
\sqrt{\sum_{\ell=1}^n (\rho_\ell-\rho)\cdot \rho_\ell^2\cdot (\xi^{(\ell)}-\eta{(\ell)})^2}
\cdot
\sqrt{\sum_{k=1}^n
(\rho_k-\rho)\cdot |\nabla^{(k)}\log h(y)|^2}.
\endalign
$$
It follows that for all $y$, $\xi$ and $\eta$
$$
\align
&\biggl|\sum_{k,\ell\in [1,n], \ell\neq k}
\rho_\ell\cdot
\biggl[\nabla^{(k)}V\bigl(\bar y^{(\ell)}, \xi^{(\ell)}\bigr)
- \nabla^{(k)} V\bigl(\bar y^{(\ell)},\eta^{(\ell)}\bigr)\biggr]
\cdot \nabla^{(k)} \log h(y)\biggr|\\
&=
\biggl|U(y,\xi)-U(y,\eta)\biggr|\\
&\leq
\sqrt{\sum_{\ell=1}^n
(\rho_\ell-\rho)\cdot \rho_\ell^2\cdot \bigl|\xi^{(\ell)}-\eta^{(\ell)}\bigr|^2}
\cdot
\sqrt{\sum_{k=1}^n
(\rho_k-\rho)\cdot\bigl|\nabla^{(k)} \log h(y)\bigr|^2}.
\endalign
$$
\bsk

Now the last line of (3.18) can be estimated as follows:

$$
\align
&\Bbb E_{p,\Pi}
\biggl|\biggl\{
\sum_{k,\ell\in [1,n], \ell\neq k}
\rho_\ell\cdot
\bigl[ \nabla^{(k)}V\bigl(\bar y^{(\ell)}, \xi^{(\ell)}\bigr)
 -\nabla^{(k)}V\bigl(\bar y^{(\ell)},\eta^{(\ell)}\bigr)\bigr]
\cdot  \nabla^{(k)}  \log h
\biggr\}\biggr|\\
&\leq
\sqrt{\sum_{\ell=1}^n
(\rho_\ell-\rho)\cdot \rho_\ell^2\cdot
\Bbb E_{p,\Pi} \bigl|\xi^{(\ell)}-\eta^{(\ell)}\bigr|^2}
\cdot
\sqrt{\sum_{k=1}^n
(\rho_k-\rho)\cdot \Bbb E_p\bigl|\nabla^{(k)} \log h(y)\bigr|^2}.
\tag 3.20
\endalign
$$
\bsk

Our calculations are valid for any coupling of
the conditional densities
$p^{(\ell)}(d\xi^{(\ell)}|\bar y^{(\ell)})$ and $Q^{(\ell)}(d\eta^{(\ell)}|\bar y^{(\ell)})$.
Now we specify $\pi^{(\ell)}(d\xi^{(\ell)},d\eta^{(\ell)}|\bar y^{(\ell)})$
so as to achieve
$$
\Bbb E_{\pi_\ell}\bigl\{ |\eta^{(\ell)}-\xi^{(\ell)}|^2|\bar y^{(\ell)}\bigr\}
=
W^2\biggl( p^{(\ell)}(\cdot|\bar y^{(\ell)}),
Q^{(\ell)}(\cdot|\bar y^{(\ell)})\biggr)\qute{for any}\quad \bar y^{(\ell)}.
$$
\bsk

By Assumptions 1 and 2, the Otto-Villani theorem can be applied to $Q^{(\ell)}(\cdot|\bar y^{(\ell)})$.  Using also the \LSI\
for $Q^{(\ell)}(\cdot|\bar y^{(\ell)})$, we get
$$
\align
&\Bbb E_{\pi_\ell}\bigl\{ |\eta^{(\ell)}-\xi^{(\ell)}|^2\bigm|\bar y^{(\ell)}\bigr\}
\leq
\frac 2{\rho_\ell}\cdot D\bigl(p^{(\ell)}(\cdot|\bar y^{(\ell)})
\bV Q_i(\cdot|\bar y^{(\ell)})\bigr)\\
&\leq
\frac 1 {\rho_\ell^2}\cdot I\biggl(p^{(\ell)}\bigl(\cdot|\bar y^{(\ell)}\bigr)
\bV Q^{(\ell)}\bigl(\cdot|\bar y^{(\ell)}\bigr)\biggr)
=
\frac 1 {\rho_\ell^2}\cdot\Bbb E_p \biggl\{\bigl|\nabla^{(\ell)} \log h\bigr|^2
  \bigm|\bar y^{(\ell)}\biggr\}
\tag 3.21
\endalign
$$
for any $\bar y^{(\ell)}$.
Substituting (3.21) into (3.20):
$$
\align
&\Bbb E_{p,\Pi}
\biggl|\biggl\{
\sum_{k,\ell\in [1,n] \ell\neq k}
\rho_\ell\cdot
\bigl[ \nabla^{(k)}V\bigl(\bar y^{(\ell)}, \xi^{(\ell)}\bigr)
 -\nabla^{(k)}V\bigl(\bar y^{(\ell)},\eta^{(\ell)}\bigr)\bigr]
\cdot  \nabla^{(k)}  \log h
\biggr\}\biggr|\\
&\leq
\sumk (\rho_k-\rho)\cdot \Bbb E_p \bigl|\nabla^{(k)} \log h\bigr|^2.
\tag 3.22
\endalign
$$
Substituting (3.19) and  (3.22) into (3.18):
$$
\align
&\Bbb E_p\bigl\{\nabla \chi\cdot \nabla \log h\bigr\}
\ge
\sum_{k=1}^n \rho_k\cdot \Bbb E_p \bigl|\nabla^{(k)} \log h\bigr|^2
-
\sum_{k=1}^n
(\rho_k-\rho)
\cdot
\Bbb E_p \bigl|\nabla^{(k)} \log h\bigr|^2\\
&=
\rho\cdot
\sum_{k=1}^n
\bigl|\nabla^{(k)} \log h\bigr|^2
=
\rho\cdot
\Bbb E_p \bigl|\nabla\log h\bigr|^2.
\qquad\qquad\qquad\qquad\qquad\qquad\qquad\qquad\qed
\endalign
$$

\enddemo

\beginsection 4. {Proof of the Auxiliary Theorem and the Approximation Lemma}

In the proof of the Auxiliary Theorem we use  the weighted
Gibbs sampler $\Ga$ with weights $\rho_k/R$, $R=\sum_k \rho_k$,
defined in Section 2:
$$
\Ga=\frac 1 R\cdot\sumk \rho_k\cdot \Ga_k,
\qquad
\Gamma_k(z|y) =\de\bigl(\bar y^{(k)}, \bar z^{(k)}\bigr)\cdot
Q^{(k)}\bigl(z^{(k)}|\bar y^{(k)}\bigr).
$$
 Recall that
 $$
\sup_{x,\xi}\bigl\Vert A(x,\xi) \bigr\Vert\triangleq 1-\de<1.
$$
\bsk

\proclaim{Proposition 4}
Under Assumptions  1-3,
for fixed $z,u\in\R^N$ we have
$$
\align
&\sum_{k=1}^n
\rho_k \cdot W^2\biggl(Q^{(k)}\bigl(\cdot|\bar z^{(k)}\bigr),
 Q^{(k)}\bigl(\cdot|\bar u^{(k)}\bigr)\biggr)
\leq
2\cdot
\sum_{k=1}^n
 D\biggl(Q^{(k)}\bigl(\cdot|\bar z^{(k)}\bigr)
 \bV Q^{(k)}\bigl(\cdot|\bar u^{(k)}\bigr)\biggr)\\
&\leq
(1-\de)^2
\cdot
\sum_{k=1}^n \rho_k\cdot  \bigl|z^{(k)}- u^{(k)}\bigr|^2.
\tag 4.1
\endalign
$$
\endproclaim

\demo{Proof}
\newline
The first \ineq\ follows from the
 Otto-Villani theorem  for $Q^{(k)}\bigl(\cdot|\bar u^{(k)}\bigr)$.
Then we use  the \LSI\
to continue (4.1) as follows:
$$
\align
&\leq
\sum_{k=1}^n \frac 1 {\rho_k}\cdot
I\biggl(
Q^{(k)}\bigl (\cdot|\bar z^{(k)}\bigr)
\bV  Q^{(k)}\bigl (\cdot|\bar u^{(k)}\bigr)
\biggr)\\
&=
\int_{\R^N}
\sum_{k=1}^n
\frac 1 {\rho_k}\cdot
\biggl|
\nabla^{(k)} V\bigl(\bar z^{(k)},\eta^{(k)}\bigr)
-
\nabla^{(k)} V\bigl(\bar u^{(k)},\eta^{(k)}\bigr)
\biggr|^2
\prod_{i=1}^n Q^{(k)}(d\eta^{(k)}|\bar z^{(k)}).
\tag 4.2
\endalign
$$
\bsk

To estimate the sum under the integral in  (4.2), fix $\eta^N$, and
consider the  function
$F=(F_1,\dots,F_N):\R^N\mapsto \R^N$
defined by
$$
\align
&F^{(k)}: \R^N\mapsto \R^{(k)},\\
&F^{(k)}(z)
=
\frac 1 {\sqrt{\rho_k}}\cdot
\nabla^{(k)} V\biggl(
\frac {z^{(1)}} {\sqrt{\rho_1}},\dots,\frac {z^{(k-1)}} {\sqrt{\rho_{k-1}}},
\frac {\eta^{(k)}} {\sqrt{\rho_k}},\frac {z^{(k+1)}} {\sqrt{\rho_{k+1}}},
\dots,\frac {z^{(n)}} {\sqrt{\rho_n}}
\biggr).
\endalign
$$
With the notation
$$
 \zeta^{(k)}=z^{(k)}\cdot \sqrt{\rho_k},\quad
\theta^{(k)}=u^{(k)}\cdot \sqrt{\rho_k},
\quad k\in [1,n],
$$
the sum under the integral in  (4.2)
is just the squared increment of $F$ between points $\zeta$ and $\theta$:
$$
\sum_{k=1}^n
\frac 1 {\rho_k}\cdot
\biggl|
\nabla^{(k)}V\bigl(\bar z^{(k)},\eta^{(k)}\bigr)
-
\nabla^{(k)}V\bigl(\bar u^{(k)},\eta^{(k)}\bigr)
\biggr|^2
=
\sumk
\biggl|F^{(k)}(\zeta)-F^{(k)}(\theta)\biggr|^2.
\tag 4.3
$$
\bsk

The Jacobian of $F$ is
$$
\biggl(\frac 1 {\sqrt{\rho_k}\sqrt{\rho_\ell}}\cdot
   V_{i,j}\bigl(\bar z^{(k)},\eta^{(k)}\bigr)\biggr)_{i\in I_k,j\in I_\ell,k\neq \ell}.
$$
(It has zeros for $i$ and $j$ belonging to the same $I_k$.)
Thus, by (2.4),
$$
\sum_{k=1}^n \bigl|F_k(\zeta)-F_k(\theta)\bigr|^2
\leq
(1-\de)^2
\cdot
\sum_{k=1}^n \rho_k\cdot \bigl|z^{(k)}- u^{(k)}\bigr|^2.
\tag 4.4
$$
Substituting (4.3) and (4.4) into (4.2) we get the desired result (4.1).
$\qquad\qquad\qed$
\enddemo
\bsk

We use Proposition 4 to show that the Gibbs sampler $\Gamma$
 is a contraction \wre\
a weighted  Wasserstein distance.
\bsk

\definition{Definition}
Let $r$ and  $s$ probability measures $r$ and  $s$ on  $\R^N$.
We define the  weighted quadratic Wasserstein distance
of $r$ and  $s$ (with wights $\rho_k$)
by
$$
W^2_{\{\rho_k\}}(r,s)  = \inf_\pi
\sum_{k=1}^n \rho_k \cdot
\Bbb E_\pi \bigl|Z^{(k)}-U^{(k)}\bigr|^2,
$$
where $Z$ and $U$ are random \sqe s s with laws $r$ resp. $s$, and infimum
is taken over all distributions $\pi=\Lw(Z,U)$  with marginals $r$
and $s$.
\enddefinition
\bsk

\proclaim{Proposition 5}
\newline
If Assumptions 1-3 hold for $q$ then
$$
W_{\{\rho_k\}}\bigl(r\Ga,s\Ga)
\leq
\biggl(1-\frac {\rho_{min}\cdot\de}R \biggr)
\cdot
W_{\{\rho_k\}}\bigl(r,s\bigr).
\tag 4.5
$$
\endproclaim
\bsk

\demo{Proof}
\newline
Let $Z=(Z_1,Z_2,\dots,Z_N)$ and $U=(U_1,U_2,\dots,U_N)$ be random \sqe s in $\R^N$, with
$\Lw(Z)=r$, $\Lw(U)=s$, and let
$\pi=\Lw(Z,U)$ be that joining of $r$ and $s$
that achieves $W^2_{\{\rho_k\}}\bigl(r,s\bigr)$.
Select a random index $\kappa\in[1,n]$ according to the distribution $(\rho_k/R)$,
and define
$$
\Lw(Z'|Z,U)=\Ga_\kappa(\cdot|Z),\quad \Lw(U'|Z,U)=\Ga_\kappa(\cdot|U).
$$
Then $\Lw(Z')=r\Ga$, and $\Lw(U')=s\Ga$.
Further, define $\Lw(Z',U')$
 as that coupling of $r\Ga$ and $s\Ga$
that achieves $W^2\bigl(Q^{(\kappa)}(\cdot|\bar Z^{(\kappa)}),
Q^{(\kappa)}(\cdot|\bar U^{(\kappa)})\bigr)$
for each value of the condition. Thereby we have defined
$\Lw(Z',U'|Z,U)$, and by Proposition 4 we have
$$
\align
&W^2_{\{\rho_k\}}\bigl(r\Ga,s\Ga)
\leq
\sum_{k=1}^n \rho_k\cdot
\biggl(1-\frac {\rho_k} R+\frac {\rho_k} R\cdot (1-\de)^2\biggr)
\cdot
\Bbb E\bigl|{Z}^{(k)}-{U}^{(k)}\bigr|^2\\
&=
\biggl(1-\frac {2\rho_{min}\cdot\de\cdot(1-\de/2)} R\biggr)
\cdot \sum_{k=1}^n \rho_k\cdot \Bbb E|Z^{(k)}-U^{(k)}\bigr|^2\\
&\leq
\biggl(1-\frac {\rho_{min}\cdot\de}{R}\biggr)
\cdot
W^2_{\{\rho_k\}}\bigl(r,s).
\qquad\qquad\qquad\qquad\qquad\qquad\qquad\qquad\qquad\qed
\endalign
$$
\enddemo
\bsk

In the sequel we shall use the
\bsk

\definition{Notation}
$$
I\bigl(p^{(k)}(\cdot|\bar Y^{(k)}) \bV Q^{(k)}(\cdot|\bar Y^{(k)})\bigr)
\triangleq
\Bbb E
I\bigl(p^{(k)}(\cdot|\bar Y^{(k)}) \bV Q_i(\cdot|\bar Y^{(k)})\bigr)
$$
(omitting the symbol of expectation).
\enddefinition
\bsk

\proclaim{Proposition 6}
\newline
Under Assumptions  1-3 we have
$$
\align
&W_{\{\rho_k\}}(p,q)
\leq
\frac{2R} {\rho_{min}\cdot\de}
\cdot
\sqrt{
\sum_{k=1}^n \rho_k \cdot
\Bbb E W^2\biggl(p^{(k)}\bigl(\cdot|\bar Y^{(k)}\bigr),
Q^{(k)}\bigl(\cdot|\bar Y^{(k)}\bigr)\biggr)}\\
&\leq
\frac  {2R}{\rho_{min}\cdot\de}
\cdot
\sqrt{
\sum_{k=1}^n
\frac 1 \rho_k \cdot
I\biggl(p^{(k)}\bigl(\cdot|\bar Y^{(k)}\bigr)\bV Q^{(k)}\bigl(\cdot|\bar Y^{(k)}\bigr)\biggr)}.
\endalign
$$
\endproclaim
\bsk

\demo{Proof}
\newline
The first inequality follows from  the triangle \ineq\ for $W_{\{\rho_k\}}(p,q)$
and Proposition 5,
and the second one follows
 from the  the Otto-Villani theorem
 and the \LSI\ for $Q^{(k)}(\cdot|\bar y^{(k)})$. $\qquad\qquad\qquad\qquad\qquad\qquad\qquad\qquad
\qed$

\enddemo
\bsk

\proclaim{Proposition 7}
\newline
There exists a $C=C(R,\rho_{min},\de)>0$ ($R=\sum_k\rho_k$ and
$\rho_{min}=\min_k\rho_k$) such that
$$
\sum_{k=1}^n  D\bigl(Y^{(k)}\bV X^{(k)}\bigr)
\leq
\frac 1 {2C}\cdot I(p||q).
\tag 4.6
$$
\endproclaim
\bsk

\demo{Proof}
\newline
Let $\pi=\Lw(Y,X)$ denote that joining of $p=\Lw(Y)$
and $q=\Lw(X)$  that achieves $W_{\{\rho_k\}}(p,q)$.
\bsk

The convexity of the entropy functional implies the \ineq\
$$
\sum_{k=1}^n  D\bigl(Y^{(k)}\bV X^{(k)}\bigr)
\leq
\sum_{k=1}^n
\Bbb E_\pi D\biggl(Y^{(k)}|\bar Y^{(k)}\bV Q^{(k)}\bigl(\cdot|\bar X^{(k)}\bigr)\biggr).
\tag 4.7
$$
\bsk

The right-hand-side of (4.7) can be written as a sum of three terms:
$$
\align
&\sum_{k=1}^n \Bbb E_\pi
D\biggl(Y^{(k)}|\bar Y^{(k)}\bV Q^{(k)}\bigl(\cdot|\bar X^{(k)}\bigr)\biggr)\\
&=
\sum_{k=1}^n  D\biggl(Y^{(k)}|\bar Y^{(k)}\bV Q^{(k)}\bigl(\cdot|\bar Y^{(k)}\bigr)\biggr)
+
\sum_{k=1}^n
\Bbb E_\pi
D\biggl(Q^{(k)}\bigl(\cdot|\bar Y^{(k)})\bV Q^{(k)}\bigl(\cdot|\bar X^{(k)}\bigr)\biggr)\\
&+
\int_{\R^N\times \R^N}
\sum_{k=1}^n
\biggl[p^{(k)}\bigl(y^{(k)}|\bar y^{(k)}\bigr)-Q^{(k)}\bigl(y^{(k)}|\bar y^{(k)}\bigr)\biggr]
\cdot
\log\frac {Q^{(k)}\bigl(y^{(k)}|\bar y^{(k)}\bigr)}
 {Q^{(k)}\bigl(y^{(k)}|\bar x^{(k)}\bigr)} dydx\\
&\triangleq S_1+S_2+S_3.
\tag 4.8
\endalign
$$
\bsk

By  the \LSI\ for
$Q^{(k)}(\cdot|\bar y^{(k)})$ we have
$$
S_1
\leq
\frac 1 2\cdot\sum_{k=1}^n
\frac 1 {\rho_k}
I\biggl(
Y^{(k)}|\bar Y^{(k)}\bV Q^{(k)}\bigl(\cdot|\bar Y^{(k)}\bigr)
\biggr).
\tag 4.9
$$
\bsk

Further, by Propositions 4,
$$
S_2
\leq
\frac{ (1-\de)^2} 2\cdot
\sum_{k=1}^n
\rho_k\cdot\Bbb E_\pi\bigl|Y^{(k)}-X^{(k)}\bigr|^2
=
\frac{ (1-\de)^2} 2\cdot W_{\{\rho_k\}}(p,q).
\tag 4.10
$$
\bsk

$S_3$ can be written as
$$
S_3
=
\Bbb E_\mu
\sum_{k=1}^n
\biggl[ V(Y)-V(\bar X^{(k)},Y^{(k)})-V(\bar Y^{(k)},\xi^{(k)})+V(\bar X^{(k)},\xi^{(k)})\biggr],
\tag 4.11
$$
where $\mu=\Lw( Y,X,\xi)$ is defined by $\Lw( Y,X)=\pi$,
$\Lw(\xi|Y,X)=\prod_{k=1}^n \Lw(\xi^{(k)}|Y)$,
and $\Lw(Y^{(k)},\xi^{(k)}|\bar Y^{(k)})$
is an arbitrary joining of $p^{(k)}(\cdot|\bar Y^{(k)})$ and $Q^{(k)}(\cdot|\bar Y^{(k)})$.
\bsk

We claim that for any quadruple of  \sqe s $(y^N,\eta^N,x^N,\xi^N)$
the following \ineq\ holds:
$$
\align
&\sum_{k=1}^n
\biggl[ V(\bar y^{(k)},\eta^{(k)})-V(\bar x^{(k)},\eta^{(k)})
-
V(\bar y^{(k)},\xi^{(k)})+V(\bar x^{(k)},\xi^{(k)})\biggr]\\
&\leq
(1-\de)\cdot
\sqrt{\sum_{k=1}^n \rho_k \bigl|y^{(k)}-x^{(k)}\bigr|^2}
\cdot
\sqrt{\sum_{k=1}^n \rho_k \bigl|\eta^{(k)}-\xi^{(k)}\bigr|^2}.
\tag 4.12
\endalign
$$
\bsk

Indeed,  introducing the function
$$
F: \R^N\times \R^N\mapsto \R,\qquad
F(y,\eta)
=
\sum_{k=1}^n
V(\bar y^{(k)},\eta^{(k)}),
$$
the left-hand-side of (4.12) can be rewritten az follows:
$$
\align
&\sum_{k=1}^n
\biggl[ V(\bar y^{(k)},\eta^{(k)})-V(\bar x^{(k)},\eta^{(k)})
-
V(\bar y^{(k)},\xi^{(k)})+V(\bar x^{(k)},\xi^{(k)})\biggr]\\
&=
F(y,\eta)-F(x,\eta)-F(y,\xi)+F(x,\xi).
\tag 4.13
\endalign
$$
\bsk

To estimate  the right-hand-side of (4.13) (with $y,x,\eta,\xi\in\R^N$  fixed), define
$$
\align
&G: [0,1]\times [0,1]\mapsto \R,\\
&G(s,t)=F\bigl(x+s(y-x),\xi+t(\eta-\xi)\bigr)\\
&=
\sum_{k=1}^n
 V\bigl(\bar x^{(k)}+s(\bar y^{(k)}- \bar x^{(k)}), \xi^{(k)}+t(\eta^{(k)}-\xi^{(k)})\bigr).
\endalign
$$
Then we have
$$
\align
&F(y,\eta)-F(x,\eta)-F(y,\xi)+F(x,\xi)\\
&=
G(1,1)-G(1,0)-G(0,1)+G(0,0).
\tag 4.14
\endalign
$$
\bsk

We have by (2.4)
$$
\align
&\Biggl|\frac{\partial^2} {\partial_s\partial_t}G(s,t)\Biggr|=\\
&\Biggl|\sum_{k,\ell\in[1,n],k\neq\ell}\quad
\sum_{i\in I_k,j\in I_\ell}
(y_i-x_i)\cdot
V_{i,j} \bigl(\bar x^{(k)} +s(\bar y^{(k)}- \bar x^{(k)}),
\xi^{(k)}+t(\eta^{(k)}-\xi^{(k)})\bigr)\cdot(\eta_j-\xi_j)\Biggr|\\
&\leq
(1-\de)\cdot
\sqrt{\sum_{k=1}^n \rho_k\cdot \bigl|y^{(k)}-x^{(k)}\bigr|^2}\cdot
\sqrt{\sum_{\ell=1}^n \rho_\ell\cdot \bigl|\eta^{(\ell)}-\xi^{(\ell)}\bigr|^2}.
\tag 4.15
\endalign
$$
Putting together (4.13), (4.14) and (4.15)  yields (4.12).
\bsk

Applying (4.12) for $\eta=y$:
$$
\align
&\sum_{k=1}^n
\biggl[ V(y)-V(\bar x^{(k)},y^{(k)})-V(\bar y^{(k)},\xi^{(k)})
+V(\bar x^{(k)},\xi^{(k)})\biggr]\\
&\leq
(1-\de)\cdot
\sqrt{\sum_{k=1}^n \rho_k \bigl|y^{(k)}-x^{(k)}\bigr|^2}\cdot
\sqrt{\sum_{\ell=1}^n \rho_\ell \bigl|y^{(\ell)}-\xi^{(\ell)}\bigr|^2}.
\tag 4.16
\endalign
$$
\bsk

Substituting (4.16) into (4.11), and using Jensen's \ineq, we get
$$
\align
&S_3
\leq
\sqrt{
\sum_{k=1}^n
\rho_k\cdot
\Bbb E \bigl|Y^{(k)}-X^{(k)}\bigr|^2}
\cdot
\sqrt{
\sum_{k=1}^n
\rho_k\cdot
\Bbb E \bigl|Y^{(k)}-\xi^{(k)}\bigr|^2}\\
&=
W_{\{\rho_k\}}(p,q)   
\cdot
\sqrt{
\sum_{k=1}^n
\rho_k\cdot
\Bbb E \bigl|Y^{(k)}-\xi^{(k)}\bigr|^2}
\tag 4.17
\endalign
$$
\bsk

To estimate the second factor, we select for $\Lw(Y^{(k)},\xi^{(k)}|\bar Y^{(k)})$
that joining of the marginals that achieves
$W^2(p^{(k)}(\cdot|\bar Y^{(k)}),Q^{(k)}(\cdot|\bar Y^{(k)}))$ for every value of the
conditions.
Then the Otto-Villani theorem
and the \LSI\ for
$Q^{(k)}(\cdot|\bar y^{(k)})$ imply the following bound for $S_3$:
$$
S_3
\leq
W_{\{\rho_k\}}(p,q)
\cdot
\sqrt{
\sum_{k=1}^n
\frac 1 \rho_k
\cdot
I\biggl(p^{(k)}\bigl(\cdot|\bar Y^{(k)}\bigr)||Q^{(k)}\bigl(\cdot|\bar Y^{(k)}\bigr)\biggr)}.
\tag 4.18
$$
\bsk

Putting together (4.9), (4.10) and (4.18):
$$
\align
&S_1+S_2+S_3\\
&\leq
\frac 1 2\cdot
\Biggl[W_{\{\rho_k\}}(p,q)
+
\sqrt{
\sum_{k=1}^n
\frac 1 \rho_k
\cdot
I\biggl(p^{(k)}\bigl(\cdot|\bar Y^{(k)}\bigr)||Q^{(k)}\bigl(\cdot|\bar Y^{(k)}\bigr)\biggr)}
\Biggr]^2.
\tag 4.19
\endalign
$$
(4.7), together with (4.8) and (4.19), completes the proof of Proposition 7.
 $\qquad\qed$

\enddemo
\bsk

\demo{Proof of the Auxiliary Theorem }
\newline
The proof goes by induction on $n$.
It is clear that for any $k$ and  $y^{(k)}\in \R^{(k)}$,
Assumptions 1-3 formulated before  Theorem 1 do hold for $n=1$, $N=|I_k|$ and
the distribution $Q^{(k)}(\cdot|y^{(k)})$.
Assume that we have proved the  Auxiliary Theorem for $n-1$ in place of $n$.
\bsk

By a well known identity for relative entropy, we have
$$
\align
&D(p||q)
=
D\bigl(Y||X\bigr)\\
&=
\frac 1 n
\sum_{k=1}^n  D\bigl(Y^{(k)}\bV X^{(k)}\bigr)
+
\frac 1 n
\sum_{k=1}^n
D\biggl(\bar Y^{(k)}|Y^{(k)}\bV \bar q^{(k)}\bigl(\cdot|Y^{(k)}\bigr)\biggr).
\tag 4.20
\endalign
$$
\bsk

Assume  the Auxiliary Theorem  for $n-1$.
By the induction hypothesis,
$$
\align
&D\biggl(\bar Y^{(k)}|Y^{(k)}\bV \bar q^{(k)}\bigl(\cdot|Y^{(k)}\bigr)\biggr)
\leq
\frac 1 {2C}
\cdot
\sum_{\ell\neq k}
I\biggl(Y^{(\ell)}|\bar Y^{(\ell)}\bV Q^{(\ell)}\bigl(\cdot|Y^{(\ell)}\bigr)\biggr)
\qute{for all}\quad k.
\endalign
$$
Thus
$$
\align
&\frac 1 n
\sum_{k=1}^n
D\bigl(\bar Y^{(k)}|Y^{(k)}\bV \bar q^{(k)}(\cdot|Y^{(k)})\bigr)\\
&\leq
\bigl(1-1/n\bigr)\cdot \frac 1 {2C}
\sum_{\ell=1}^n
I\biggl(Y^{(\ell)}|\bar Y^{(\ell)}\bV Q^{(\ell)}\bigl(\cdot|Y^{(\ell)}\bigr)\biggr)\\
&=
\bigl(1-1/n\bigr)\cdot \frac 1 {2C}\cdot I(p||q).
\tag 4.21
\endalign
$$
\bsk

Substituting (4.6) (Proposition 7) and (4.21) into (4.20)
 completes the proof of the Auxiliary Theorem.
$\qquad\qquad\qquad\qquad\qquad\qquad\qquad\qquad\qquad\qquad\qquad
\qed$
\enddemo

\demo{Proof of the Approximation Lemma}
\newline
First we keep $q$ fixed, and
construct
 a density $g\in\Cal C^\infty$ with compact support, and such that,
with the notations
$$
r=g\cdot q, \quad
\bar r^{(k)}=\int_{\R^{(k)}} r(\bar y^{(k)},\xi^{(k)})d\xi^{(k)},
$$
we have
\bsk

$$
 D(r||q)\qute{is arbitrarily close to}\quad D(p||q),
\tag 4.22
$$
and
$$
D\bigl(\bar r^{(k)} ||\bar q^{(k)}\bigr)
\qute{is arbitrarily close to}\quad D(\bar p^{(k)}||\bar q^{(k)}\bigr)
\qute{for all}\quad k\in[1,n].
\tag 4.23
$$
\bsk

Denote by $B_m$ the closed ball in $\R^N$ around the origin and with radius $m$.
Let $\phi_m: \R^N\mapsto [0,1]$ be a $\Cal C^\infty$ function
satisfying
$$
\phi_m(x)=1\qute{for}\quad x\in B_m,
\quad
\phi_m(x)=0\qute{for}\quad x\notin B_{m+1}.
$$
Set
$$
g_m(x)=\frac 1 {\alpha_m}\cdot h(x)\cdot \phi_m(x),
\qute{and}
\quad r_m(x)=g_m(x)\cdot q(x),
$$
where $\alpha_m=\int_{\R^N} h(x)\cdot \phi_m(x)q(dx)$.

\bsk

We have
$$
D(r_m||q)
=
\frac 1 {\alpha_m}
\int_{\R^N}
\biggl(h(x) \phi_m(x)\biggr)\cdot\log\biggl(h(x) \phi_m(x)\biggr) q(dx)
-
\log \alpha_m,
$$
and $\lim_{t\to 1}\alpha_m=1$.
Since $h(x)\cdot \phi_m(x)\to h(x)$ everywhere, with
$\bigl|\bigl(h(x) \phi_m(x)\bigr)\cdot\log\bigl(h(x) \phi_m(x)\bigr)\bigr|_+$ increasing,
and using also the \ineq\  $u\log u\ge-1/e$,
 it follows that
$$
\lim_{m\to\infty}
\int_{\R^N}
\biggl(h(x) \phi_m(x)\biggr)\cdot\log\biggl(h(x) \phi_m(x)\biggr) q(dx)
=
D(p||q).
$$
Putting $g=g_m$ and $r=r_m$, for large enough $m$ we achieve (4.22).
It can be proved similarly that (4.23) can be achieved as well.
\bsk

Again, it is easily seen
that
$$
\lim_{\e\to 0}D\bigg(\bigl(1-\e)\cdot r+\e\cdot q\bigr)\bV q\biggr)
=
D(r||q),
$$
and
$$
\align
&\lim_{\e\to 0} D\bigg(\bigl(1-\e)\cdot r^{(k)}(\cdot|\bar Y^{(k)}\bigr)
+
\e\cdot q^{(k)}(\cdot|\bar Y^{(k)}\bigr)
\bV
q^{(k)}(\cdot|\bar Y^{(k)}\bigr)\biggr)\\
&=
D\bigl(r^{(k)}(\cdot|\bar Y^{(k)}\bigr)\bV q^{(k)}(\cdot|\bar Y^{(k)}\bigr).
\endalign
$$
Thus, for $q$ fixed, $h$ can be replaced  by $f=(1-\e)g+\e$.
\bsk

Now we can assume that $h$ is of the form claimed in the
 Approximation lemma. We keep the notation $p=h\cdot q$
 with the newly defined $h$, and keep  $h$ fixed.
\bsk

Now we approximate $q(x)$ by an increasing \sqe\ $\tilde q_m(x)\in\Cal C^\infty$,
and set
$$
q_m=\frac{\tilde q_m}{\int \tilde q_m(x)dx}\qute{and}\quad
\bar q_m^{(k)}(\bar x^{(k)})
=
\int_{\R^{(k)}}q_m(x)dx^{(k)}.
$$
Then define  $p_m(x)=h(x)\cdot q_m(x)$.
Since $h$ is smooth and bounded from below and above, it is easily seen that
$$
D(p_m||q_m)=\int_{\R^N}h(x)\log h(x)q_m(dx)\to D(p||q),
$$
and
$$
D\bigl(\bar p_m^{(k)}\bV\bar q_m^{(k)}\bigr)\biggr)
=
\int_{\R^{(k)}}
\bar h^{(k)}(\bar x^{(k)})\cdot
\log
\bar h^{(k)}(\bar x^{(k)})\bar q_m^{(k)}(d\bar x^{(k)})
\to
D(\bar p^{(k)}||\bar q^{(k)}).
$$
This completes the proof of the Approximation Lemma.
$\qquad\qquad\qquad\qquad\qquad\qquad\qed$

\enddemo

\beginsection{Acknowledgment}

The author is grateful for the patient support and help by P. E. Frenkel.
\bsk
\bsk
\bsk

\enddocument